\documentclass[11pt]{article}
\usepackage{amsmath}
\usepackage{amssymb}
\usepackage{enumerate}
\usepackage{color}
\usepackage{xfrac}
\usepackage{bbm}
\usepackage{hyperref}
\usepackage{mathrsfs}

\newtheorem{thm}{Theorem}[section]
\newtheorem{defn}[thm]{Definition}
\newtheorem{lemma}[thm]{Lemma}
\newtheorem{prop}[thm]{Proposition}
\newtheorem{rem}[thm]{Remark}
\newtheorem{cor}[thm]{Corollary}

\newcommand{\bl}{\mathcal B_{\lambda}}
\newcommand{\tla}{\mathcal T_{\lambda,\alpha}}

\begin{document}

\title{Some two-parameter families of generalized Bernstein functions}

\author{Stamatis Koumandos and
Henrik L. Pedersen\footnote{Research supported by grant DFF-1026-00267B from The Danish Council for Independent Research $|$ Natural Sciences}
}

\date{\today}
\maketitle

\begin{abstract}
We prove that certain functions involving ratios of Gamma functions and the Psi-function belong to generalized Bernstein classes and new properties of generalized Bernstein functions are given. 
\end{abstract}


\section{Introduction}
The behaviour of ratios of Gamma and related functions has attracted the attention of several authors over the years, see e.g.\ \cite{AP,BK, BKP,ism}. In this paper we investigate properties of the family of functions 
\begin{equation}
    \label{eq:def-mst}
    M_{a,b}(x)=x^{a-b}\frac{\Gamma(x+b)}{\Gamma(x+a)}, \;\; x>0,
\end{equation}
for positive parameters $a$ and $b$. 
This family was studied in \cite{KL}. 
We shall, in particular, relate these functions to classes of higher order Thorin-Bernstein functions and classes of generalized Stieltjes functions of positive order. We generalize results for the function $g_{\lambda}(x)=x^{\lambda}\Gamma(x)/\Gamma(x+\lambda)=M_{\lambda,1}(x)$, investigated in connection with sharp inequalities involving the incomplete Beta and Gamma functions, see \cite{KP5}.

Because of the symmetry $M_{a,b}=1/M_{b,a}$ we shall assume that
$0<b<a$
throughout the paper.

For the reader's convenience we give some fundamental definitions. A function $f:(0,\infty)\to \mathbb R$ is called completely monotonic if $f$ is infinitely often differentiable and $(-1)^nf^{(n)}(x)\geq 0$ for $n=0,1,\ldots$ and $x>0$. A $C^{\infty}$-function $f:(0,\infty)\to (0,\infty)$ belongs to the class $\bl$ of generalized Bernstein functions of order $\lambda>0$ if $x^{1-\lambda}f'(x)$ is completely monotonic. This class was investigated in \cite{KP1} where it was shown that $f$ belongs to $\mathcal B_{\lambda}$ if and only if it admits a representation 
\begin{equation}
\label{eq:bl}
f(x)=ax^{\lambda}+b+\int_0^{\infty}\gamma(\lambda,xt)\, \frac{d\mu(t)}{t^{\lambda}},
\end{equation}
where $a,b\geq 0$ and $\mu$ is a positive measure on $[0,\infty)$ making the integral converge. Here $\gamma$ denotes the incomplete Gamma function 
$$
\gamma(\lambda,x)=\int_0^{x}t^{\lambda-1}e^{-t}\, dt.
$$
Notice also that $\mu$ in \eqref{eq:bl} is the restriction to $(0,\infty)$ of the measure $\tilde{\mu}$ given in Bernstein's representation of the completely monotonic function $f'(x)x^{1-\lambda}$, i.e.\ by $\mathcal L(\tilde{\mu})(x)=f'(x)x^{1-\lambda}$, where $\mathcal L$ denotes the Laplace transform.

We remark that when $\lambda=1$ the class of functions, denoted by $\mathcal{B}$, is the so-called Bernstein functions appearing in probability theory. See \cite{S}. 

The class $\mathcal T_{\lambda,\alpha}$ (where $\lambda>0$ and $\alpha<\lambda+1$) of higher order Thorin-Bernstein functions is the subclass of $\mathcal B_{\lambda}$ defined as follows: $f\in \mathcal T_{\lambda,\alpha}$ if 
$$
f(x)=ax^{\lambda}+b+\int_0^{\infty}\gamma(\lambda,xt)\varphi(t)\, dt,
$$
where $\varphi$ is a completely monotonic function of order $\alpha$ (meaning that $t^{\alpha}\varphi(t)$ is completely monotonic) such that the integral converges.

A function $g:(0,\infty)\to \mathbb R$ belongs to the class $\mathcal S_{\lambda}$ of generalized Stieltjes functions of order $\lambda>0$ if 
$$
g(x)=\int_0^{\infty}\frac{d\mu(t)}{(x+t)^{\lambda}}+c
$$
for some positive measure $\mu$ on $[0,\infty)$ and some $c\geq 0$. An equivalent formulation is that $g$ admits the representation
\begin{equation}
\label{eq:alt-stieltjes}
g(x)=\int_0^{\infty}e^{-xt}t^{\lambda-1}\varphi(t)\, dt+c,
\end{equation}
where $\varphi$ is completely monotonic. Let us also mention that products of generalized Stieltjes functions are again generalized Stieltjes functions in the sense that 
$$
\mathcal S_{\lambda_1}\mathcal S_{\lambda_2}\subset \mathcal S_{\lambda_1+\lambda_2},
$$
see \cite{HW}.

There is a close relation between the higher order Thorin-Bernstein functions and the generalized Stieltjes functions of positive order. Indeed, according to  \cite[Theorem 2.14]{KP4} we have
 $$ f\in \tla \ \Leftrightarrow \ x^{1-\lambda}f'(x)\in \mathcal S_{\lambda+1-\alpha}.
 $$

 We give conditions on the parameters $a$ and $b$ so that the function $M_{a,b}$ belongs to $\mathcal T_{a-b,-3}$ and in addition, the function $L_{a,b}(x)=x(1-M_{a,b}(x))$ is shown to be a generalized Bernstein function of order $a-b$ for $a-b>1$. In the final section we present new properties on generalized Bernstein functions.

\section{Preliminary results}
Our results are based on integral representations of $\log M_{a,b}$ and its derivative in terms of the function function $\Phi_{a,b}$ defined as
$$
\Phi_{a,b}(t)=(b-a)+\frac{e^{-bt}-e^{-at}}{1-e^{-t}}.
$$
We begin by giving representations of $\Phi_{a,b}$ and its derivative. 

The point mass at the point $c$ is denoted  by $\epsilon_c$ and the indicator function of a set $A$ by $1_A$. With this notation we may now state and prove the first lemma.

\begin{lemma}
\label{lemma:Phi} The representation
\begin{align*}
    -\frac{\Phi_{a,b}(t)}{t^2}&=\mathcal L\left(\xi\right)(t)
\end{align*}
holds with
$$
 \xi(s)=(a-b)s+\sum_{k=0}^{\infty}\left( 1_{[k+a,\infty)}(s)(s-k-a)-1_{[k+b,\infty)}(s)(s-k-b)\right).
 $$
 Furthermore, 
 \begin{align*}
    -\frac{\Phi_{a,b}'(t)}{t^2}&=\mathcal L\left(\eta\right)(t)
\end{align*}
holds with
$$
 \eta(s)=\sum_{k=0}^{\infty}\left( (b+k)1_{[k+b,\infty)}(s)(s-k-b)-(a+k)1_{[k+a,\infty)}(s)(s-k-a)\right).
 $$
\end{lemma}
{\it Proof.}
Since
$$
\frac{1}{1-e^{-t}}=\mathcal L\left(\sum_{k=0}^{\infty}\epsilon_k\right)(t)
$$
we see that 
\begin{align}
 -\Phi_{a,b}(t)&=
 \mathcal L\left((a-b)\epsilon_0+\sum_{k=0}^{\infty}(\epsilon_{k+a}-\epsilon_{k+b})\right)(t).
 \label{eq:Phi}
 \end{align}
 Therefore (using $sds\ast\epsilon_c=1_{[c,\infty)}(s)(s-c)$),
\begin{align*}
 -\frac{\Phi_{a,b}(t)}{t^2}&=\mathcal L \left(sds\right)(t)
 \mathcal L\left((a-b)\epsilon_0+\sum_{k=0}^{\infty}(\epsilon_{k+a}-\epsilon_{k+b})\right)(t)\\
 &=\mathcal L\left(sds\ast \left((a-b)\epsilon_0+\sum_{k=0}^{\infty}(\epsilon_{k+a}-\epsilon_{k+b})\right)\right)(t)
=\mathcal L \left(\xi\right)(t),
 \end{align*}
 with $\xi$ as in the statement of the lemma. 

 To prove the second assertion we differentiate in \eqref{eq:Phi} and obtain
 \begin{align*}
 -\Phi_{a,b}'(t)&=
 \mathcal L\left(\sum_{k=0}^{\infty}((k+b)\epsilon_{k+b}-(k+a)\epsilon_{k+a})\right)(t).
 \end{align*}
 Hence
  \begin{align*}
    -\frac{\Phi_{a,b}'(t)}{t^2}&=\mathcal L \left(sds\ast \sum_{k=0}^{\infty}((k+b)\epsilon_{k+b}-(k+a)\epsilon_{k+a})\right)(t)=\mathcal L(\eta)(t).
\end{align*}
This completes the proof.\hfill$\square$

\begin{figure}
\begin{center}
       \includegraphics[width=0.5\textwidth]{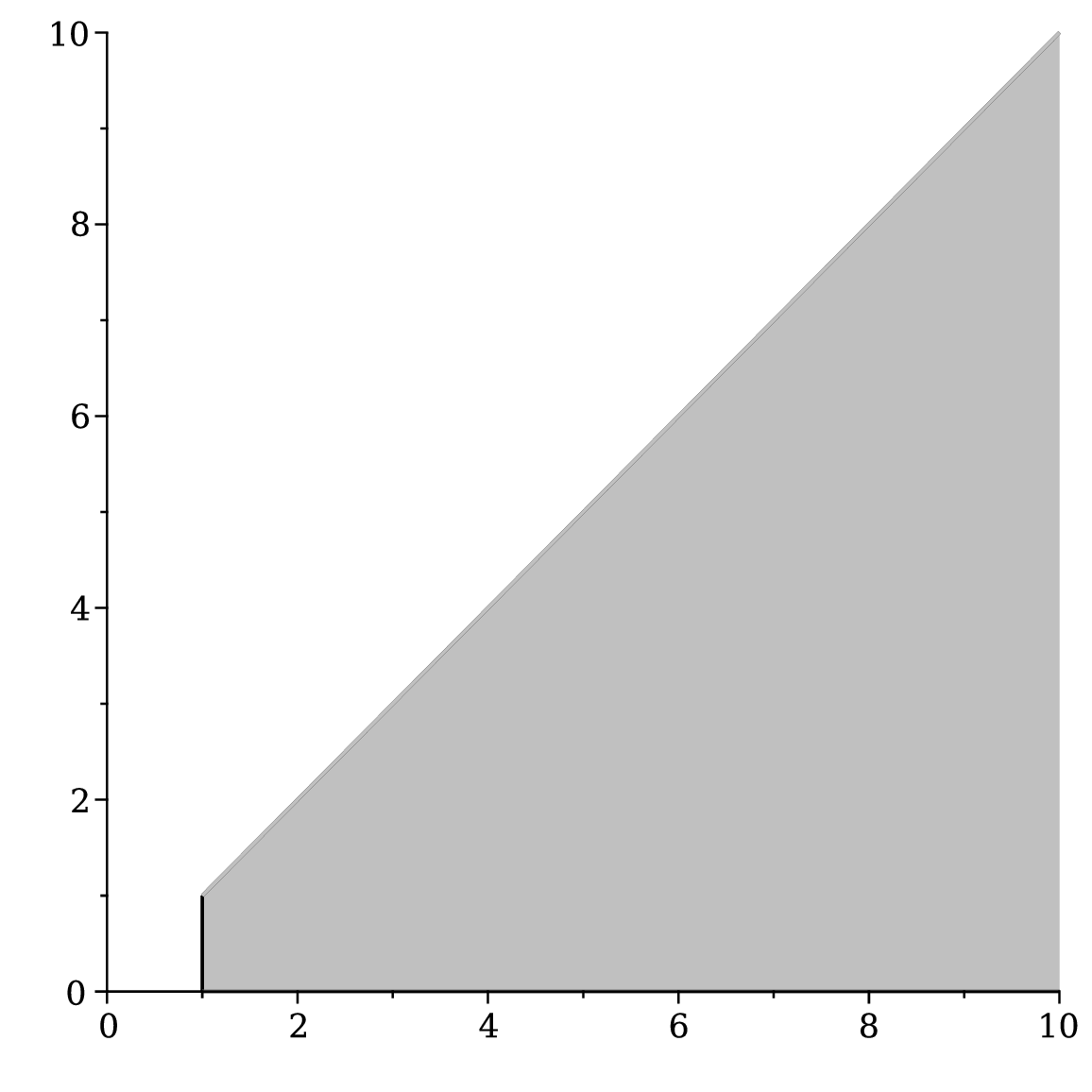}
 \end{center}  
 \caption{The parameter domain $\Omega$}
 \label{fig:omega}
\end{figure}

We define the domain $\Omega$ (sketched in Figure \ref{fig:omega}) as
$$
\Omega=\{(a,b)\in \mathbb R^2\, |\, 0<b<a,\ a>1\}.
$$

\begin{lemma}
    \label{lem:xi}
    When $(a,b)\in \Omega$, $\xi(t)\geq 0$ for $t\geq 0$.
\end{lemma}
{\it Proof}. It is easy to see that $\xi$ is a continuous and piecewise linear function. Furthermore, a computation shows that
 $$
 \xi(s+1)=\xi(s)+\Theta(s)
 $$
 with $\Theta(s)=(a-b)+1_{[a-1,\infty)}(s)(s+1-a)-1_{[b-1,\infty)}(s)(s+1-b)$. The function $\Theta$ is easily seen to be non-negative, and thus we need only to check that $\xi $ is non-negative on $[0,1]$. Since $a>1$ we have 
 $\xi(s)=(a-b)s-1_{[b,\infty)}(s)(s-b)$ for $s\in [0,1]$ which indeed is non-negative. \hfill $\square$
\begin{lemma}
    \label{lem:eta}
    When $(a,b)\in \Omega$, $\eta(t)\geq 0$ for $t\geq 0$.
\end{lemma}
{\it Proof}. 
If $a=b+k_0$ for some $k_0\in \{1,2,\ldots\}$ then $\eta$ is clearly non-negative on $[0,\infty)$. Otherwise, we choose $k_0\in \{1,2,\ldots\}$ such that
 $$
 b+k_0-1<a<b+k_0
 $$
 and let $b_0=b+k_0-1$. This gives
 \begin{align*}
 \lefteqn{\eta(s)=\sum_{k=0}^{k_0-2}(b+k)1_{[b+k,\infty)}(s)(s-b-k)+}\\
 &\sum_{k=0}^{\infty}\left( (b_0+k)1_{[k+b_0,\infty)}(s)(s-k-b_0)-(a+k)1_{[k+a,\infty)}(s)(s-k-a)\right).
 \end{align*}
 The first sum is clearly non-negative and the second infinite sum is a function of the same form as $\eta$, only this time with parameters $b_0, a$ such that $0<b_0<a<b_0+1$. 
 
 It follows that it is enough to prove $\eta\geq 0$ when the parameters $a$ and $b$ satisfy $b<a<b+1$. 

 Since $\eta$ is continuous, piecewise linear and $\eta(0)=0$,  
$\eta$ is non-negative provided $\eta(a+l)\geq$ and $\eta(b+l)\geq 0$ for all $l=0,1,\ldots$ (these being the points of discontinuity of the derivative). We get
\begin{align*}
    \eta(b+l)&=\sum_{k=0}^{l-1}\left((b+k)(b+l-b-k)-(a+k)(b+l-a-k)\right)\\
    &=(a-b)\sum_{k=0}^{l-1}(2k+a-l)= (a-b)l(a-1)\geq 0
\end{align*}
since $a>1$.
Similarly,
\begin{align*}
    \eta(a+l)&=\sum_{k=0}^{l}\left((b+k)(b+l-b-k)-(a+k)(b+l-a-k)\right)\\
    &=(a-b)\sum_{k=0}^{l}(2k+b-l)= (a-b)b(l+1)\geq 0.
\end{align*}
This completes the proof.
 \hfill$\square$

\begin{figure}
    \begin{center}
        \includegraphics[width=0.4\textwidth]{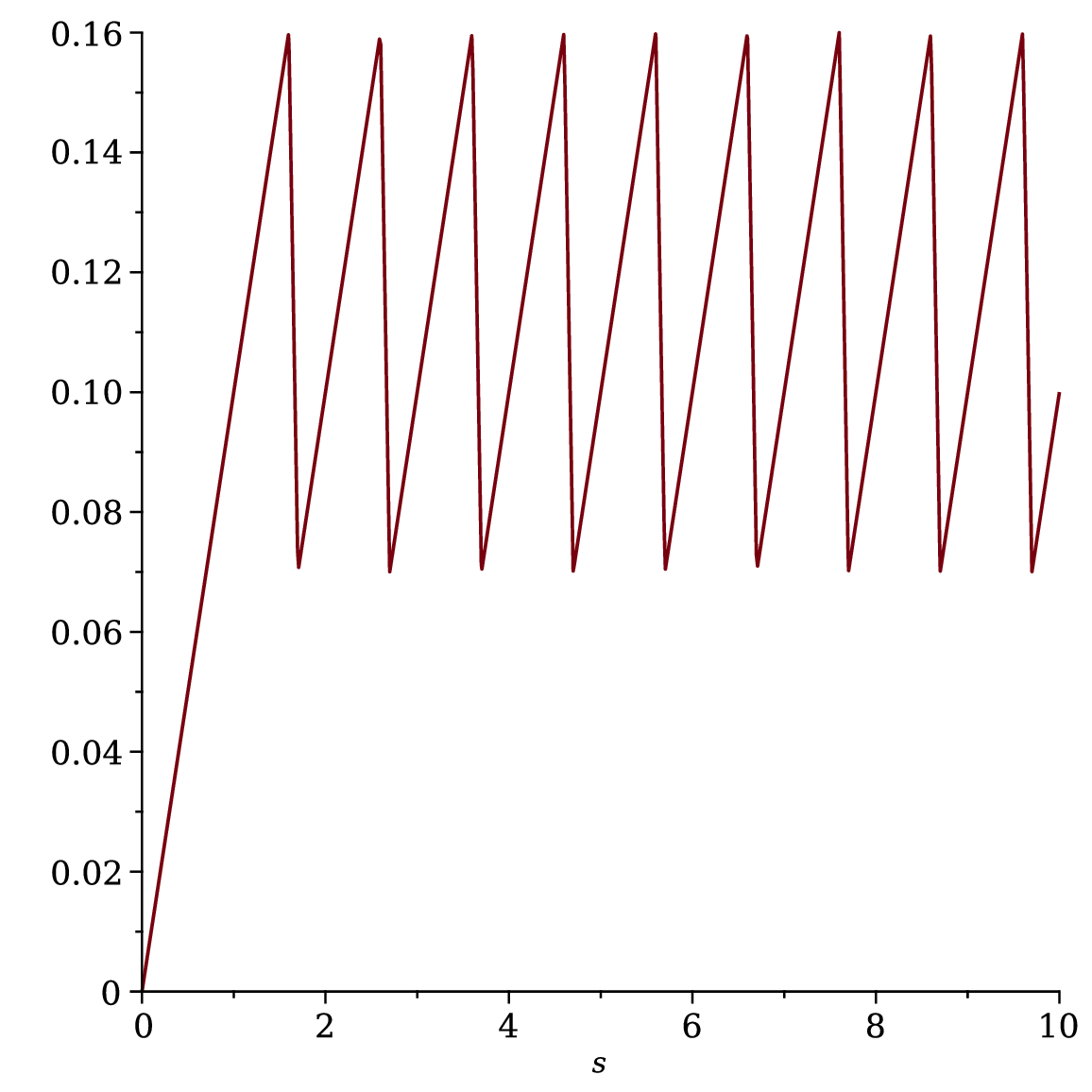}\quad
         \includegraphics[width=0.4\textwidth]{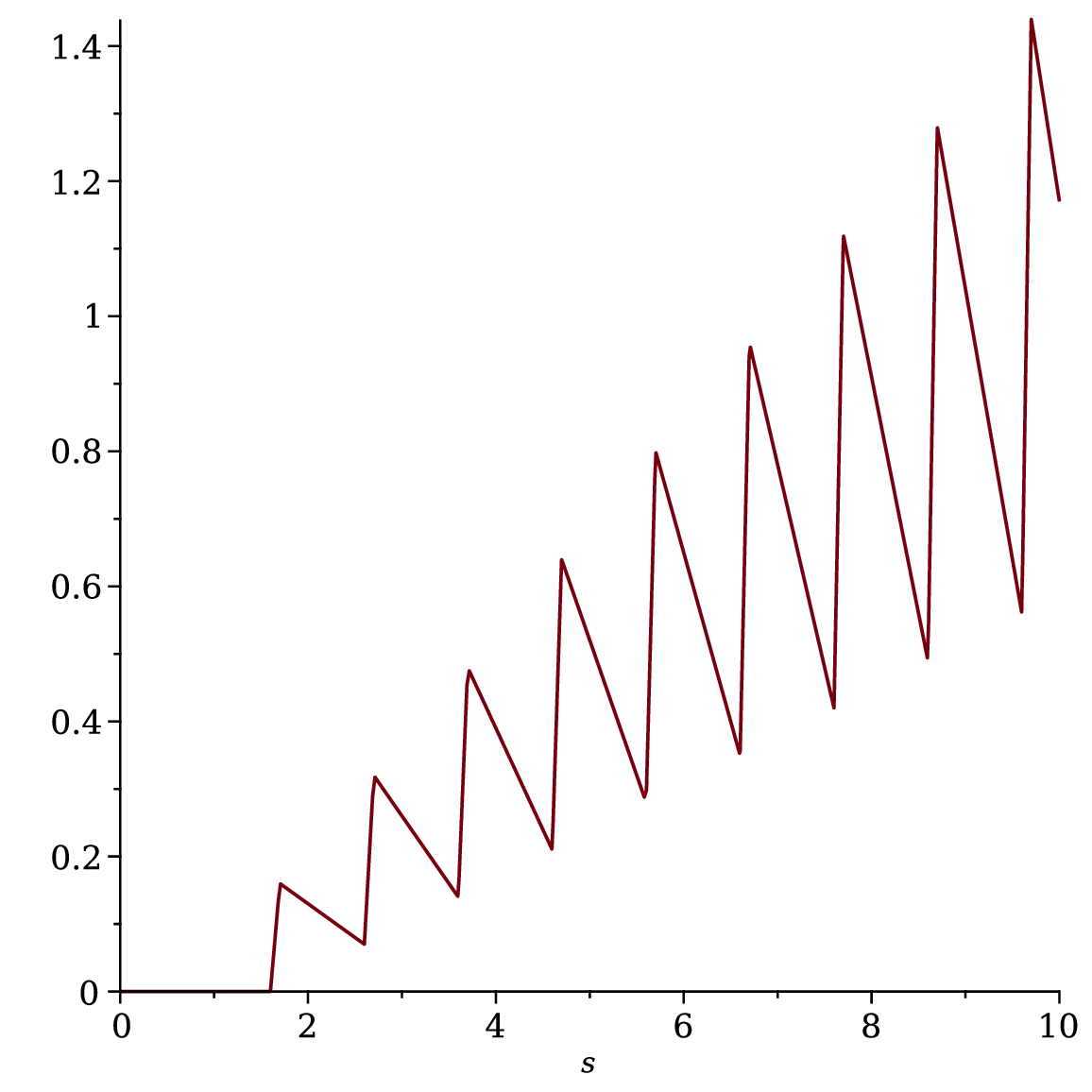}
    \end{center}
    \caption{The functions $\xi$ (to the left) and $\eta$ (to the right) corresponding to the papameters $b=1.6$ and $a=1.7$}
    \label{fig:xi+eta}
\end{figure}

 Plots of $\xi$ and $\eta$ are included in Figure \ref{fig:xi+eta}. Notice that $\xi$ is periodic for large values of the variable, since the function $\Theta$ in the proof of Lemma \ref{lem:xi} equals zero on $[a-1,\infty)$.

\section{Main results}
The integral representation of the logarithmic derivative of the Gamma function
$$
\psi(x)=(\log \Gamma)'(x)=\int_0^{\infty}\left(\frac{e^{-t}}{t}-\frac{e^{-xt}}{1-e^{-t}}\right)\, dt
$$
(see \cite{as,aar}) is used in the proof of the next lemma.
\begin{lemma}
\label{lemma:prelim}
    The following representations hold.
   \begin{align*}
     -(\log M_{a,b})'(x)&=\int_0^{\infty}e^{-xu}\Phi_{a,b}(u)\, du,\\
    -x(\log M_{a,b})'(x)&=\int_0^{\infty}e^{-xu}\Phi_{a,b}'(u)\, du,
    \\
    \log M_{a,b}(x)&=\int_0^{\infty}e^{-xu}(\Phi_{a,b}(u)/u)\, du.
   \end{align*}
\end{lemma}
{\it Proof.}
By definition we have
$$
\log M_{a,b}(x)=(a-b)\log x+\log \Gamma (x+b)-\log \Gamma(x+a),
$$
so that, using the representation of $\psi$, 
\begin{align*}
-(\log M_{a,b})'(x)&=\frac{b-a}{x}-\psi(x+b)+\psi(x+a)
=\int_0^{\infty}e^{-xt}\Phi_{a,b}(t)\, dt.
\end{align*}
Thus the first assertion holds. (This is also contained in \cite[(6.10)]{KL}, but for the reader's convenience we have included the derivation here.)  By partial integration the second assertion follows. Using Fubini's theorem and the relation $\lim_{u\to \infty}M_{a,b}(u)=1$ we obtain the third assertion 
    by integrating the first one (from $x$ to $\infty$) and this proves the lemma.\hfill$\square$

\begin{prop}
\label{prop:log-mst-in-s2}
 The function $-\log M_{a,b}$ belongs to $\mathcal S_2$ for $(a,b)\in \Omega$. \end{prop}
 {\it Proof.}
 In view of the third relation in Lemma \ref{lemma:prelim} and \eqref{eq:alt-stieltjes} we need to show that $-\Phi_{a,b}(t)/t^2$ is completely monotonic. This follows by combining Lemma \ref{lemma:Phi} and Lemma \ref{lem:xi}.
 \hfill$\square$

 \begin{rem}
 Kristiansen's theorem (see \cite{kristiansen} and \cite{BKP}) states that any function $f$ from $\mathcal S_2\setminus \{0\}$ is logarithmically completely monotonic (meaning that $-f'/f$ is completely monotonic). A function is called infinitely divisible if its Laplace transform is logarithmically completely monotonic. Combining Krisitansen's theorem with the first formula of Lemma \ref{lemma:prelim} we see that $-\Phi_{a,b}(t)/t$ is infinitely divisible. 
 \end{rem}

 \begin{prop}
     \label{prop:xlogMprime}
     The function $x(\log M_{a,b})'(x)$ belongs to $\mathcal S_3$ for $(a,b)\in \Omega$, and
     $$
     x(\log M_{a,b})'(x)=2\int_0^{\infty}\frac{\eta(t)}{(x+t)^3}\, dt.
     $$
 \end{prop}
 {\it Proof.} Using the second relation in Lemma \ref{lemma:prelim} and again \eqref{eq:alt-stieltjes} it follows that $x(\log M_{a,b})'(x)$ belongs to $\mathcal S_3$ exactly when $-\Phi_{a,b}'(t)/t^2$ is completely monotonic. This follows from Lemma \ref{lemma:Phi} and Lemma \ref{lem:eta}, and we also obtain the representation in terms of $\eta$.\hfill $\square$

 \begin{rem}
    As we have seen, the function $-\Phi_{a,b}'(t)/t^2$ is completely monotonic, and since also $-\Phi_{a,b}>0$ it follows that $-\Phi_{a,b}\in \mathcal B_{3}$. Using Lemma \ref{lemma:prelim} and \cite[Theorem 3.1]{KP1} we obtain $x(\log M_{a,b})'(x)=x\mathcal L(-\Phi_{a,b})(x)\in \mathcal S_3$. This gives another proof of Proposition \ref{prop:xlogMprime}. 
    
    Since the function $-\Phi_{a,b}$ is also bounded it has a representation 
    $$
    -\Phi_{a,b}(x)=2q(\infty)-x^3\int_0^{\infty}q(v)v^2e^{-xv}\, dv,
    $$
    where $q$ given by 
    $$
    q(t)=\int_0^t\frac{\eta(s)}{s^3}\, ds.
    $$
    is positive, increasing and bounded. (See \cite[Proposition 3.9]{KP1}.)
\end{rem}

 \begin{thm}
 \label{thm:MinT}
     For $(a,b)\in \Omega$ we have $M_{a,b}\in \mathcal T_{a-b,-3}$.
 \end{thm}
    {\it Proof.} As mentioned in the introduction it is equivalent to show that $x^{1-a+b}M_{a,b}'(x)$ belongs to $\mathcal S_{a-b+4}$. Writing $a-b=[a-b]+s$, with $s\in [0,1)$ (the square brackets denote integer part) and $t=x+b$ we obtain from  \cite[Proposition 4.6]{BKP}
$$
\frac{\Gamma(x+b)}{\Gamma(x+a)}=\frac{\Gamma(t)}{\Gamma(t+[a-b]+s)}\in \mathcal S_{[a-b]+1}.
$$
    Therefore
    \begin{align*}
      x^{1-a+b}M_{a,b}'(x)&=x^{1-a+b}M_{a,b}(x)\left(\log M_{a,b}\right)'(x)\\
      &=\frac{\Gamma(x+b)}{\Gamma(x+a)}\, x(\log M_{a,b})'(x)\in \mathcal S_{[a-b]+1}\mathcal S_3.
    \end{align*}
    Since $S_{[a-b]+1}\mathcal S_3\subseteq \mathcal S_{[a-b]+4}\subseteq\mathcal S_{a-b+4}$
    the result is proved.
     \hfill $\square$
   
Theorem \ref{thm:MinT} yields the representation 
     $$
     M_{a,b}(x)=\alpha x^{a-b}+\beta+\int_0^{\infty}\gamma(a-b,xt)t^3\varphi(t)\, dt,
     $$
     for some completely monotonic function $\varphi$. Since $M_{a,b}$ is a bounded function from $\mathcal B_{a-b}$, $\lim_{x\to \infty}M_{a,b}=1$ and $M_{a,b}(0)=0$ it thus admits the representation
     \begin{equation}
         \label{eq:Mbdd}
      M_{a,b}(x)=1-x^{a-b}\int_0^{\infty}p(t)t^{a-b-1}e^{-xt}\, dt,
     \end{equation}
     with $p(t)=\int_0^tu^3\varphi(u)\, du$ being bounded. (See \cite[Proposition 22 and 3.9]{KP1}.)

\begin{rem}
    In  \cite{KL} the family of functions $M_{a,b}$ was also investigated outside the domain $\Omega$, namely for $(a,b)$ inside the triangle with corners $(1,0)$, $(1,1)$ and $(1/2,1/2)$: for these values of the parameters we have 
    $$
    x^{1-a+b}M_{a,b}(x)=\frac{\Gamma(x+b)}{\Gamma(x+a)}x(\log M_{a,b})'(x),
    $$
    where both the Gamma-ratio and $x(\log M_{a,b})'(x)$ is completely monotonic  since $-\Phi_{a,b}'(t)\geq 0$, according to Lemma \ref{lemma:prelim}. The positivity of $-\Phi_{a,b}$ was also proved in \cite{KL}.
\end{rem}

    \section{A related function}
    In the paper \cite{KL} the function $L_{a,b}$ given by 
    $$
    L_{a,b}(x)=x(1-M_{a,b}(x))
    $$
    was investigated. 
    According to \cite[(2.10) and (2.11)]{KL} we have the concrete representation 
    \begin{equation}
    \label{eq:L}
      L_{a,b}(x)=\frac{1}{\Gamma(a-b)}x^{a-b+1}\int_0^{\infty}t^{a-b-1}(1-w(t))e^{-xt}\, dt,
    \end{equation}
    where 
    $$
    w(t)=\left(\frac{1-e^{-t}}{t}\right)^{a-b-1}e^{-bt}.
    $$
     If we compare this with \eqref{eq:Mbdd} it follows that
     $$
     \frac{1-w(t)}{\Gamma(a-b)}= \int_0^tu^3\varphi(u)\, du,
     $$
     so $-w'(t)/t^3$ is completely monotonic. (This does not seem to be obvious from the definition of $w$.)
     Furthermore we obtain the representation
     $$
     M_{a,b}(x)=-\frac{1}{\Gamma(a-b)}\int_0^{\infty}\gamma(a-b,xt)w'(t)\, dt.
     $$
    \begin{prop}
    \label{prop:L}
        For $(a,b)\in \Omega$ with $a-b>1$ we have $L_{a,b}\in \mathcal B_{a-b}\setminus \mathcal B_1$.
    \end{prop}
    {\it Proof.} From \cite[Theorem 3.1]{KL} we have $L_{a,b}\notin \mathcal B_1$ and hence we concentrate on proving $L_{a,b}\in \mathcal B_{a-b}$. 
    To ease notation let $\lambda=a-b$ and notice that $\lambda >1$.  
    Differentiation in \eqref{eq:L} and partial integration shows
    \begin{equation}
        \label{eq:Lprime}
    \Gamma(\lambda)L_{a,b}'(x)x^{1-\lambda}= x\int_0^{\infty}t^{\lambda-1}(1-w(t)+tw'(t))e^{-xt}\, dt.
    \end{equation}
    Hence $L_{a,b}'(x)x^{1-\lambda}$ is completely monotonic provided $t^{\lambda-1}(1-w(t)+tw'(t))$ is a non-negative and increasing function on $(0,\infty)$. See \cite{KP0}. In the lemma below it is proved that $w$ is  logarithmically convex. It is thus also convex and hence 
    $$
    (1-w(t)+tw'(t))'=tw''(t)>0.
    $$
    It is elementary to see that $\lim_{t\to 0}(1-w(t)+tw'(t))=0$ and therefore the function $1-w(t)+tw'(t)$ is both increasing and non-negative.\hfill$\square$

    \begin{lemma}
    The function $w$ 
    is decreasing and logarithmically convex on $(0,\infty)$ for $(a,b)\in \Omega$ with $a-b>1$.
    \end{lemma}
    {\it Proof.} Elementary calculations show that
    $$
    (\log w)''(t)=(a-b -1)\left(\frac{1}{t^2}-\frac{e^t}{(e^t-1)^2}\right).
    $$
    This expression is positive since 
    $$
    t^2<2(\cosh t-1)=e^t+e^{-t}-2=(e^t-1)(1-e^{-t})=\frac{(e^t-1)^2}{e^t},
    $$
    and the proof is complete.\hfill $\square$

\section{Additional results and comments}
We begin by a simple result about products of generalized Bernstein functions. A counterpart for higher order Thorin-Bernstein functions is proved in Corollary \ref{cor:product-higher-order}.
\begin{prop}
\label{prop:product-bernstein}
    If $f_1\in \mathcal B_{\lambda_1}$ and $f_2\in \mathcal B_{\lambda_2}$ then $f_1f_2$ belongs to $\mathcal B_{\lambda_1+\lambda_2}$.
\end{prop}
{\it Proof.} By definition,  $f_j'(x)/x^{\lambda_j-1}$ is a completely monotonic function, and so is $f_j(x)/x^{\lambda_j}$. (See \cite[Corollary 2.1]{KP1}.) The elementary relation 
$$
\frac{(f_1f_2)'(x)}{x^{\lambda_!+\lambda_2-1}}=
\frac{f_1'(x)}{x^{\lambda_1-1}}\frac{f_2(x)}{x^{\lambda_2}}+
\frac{f_1(x)}{x^{\lambda_1}}\frac{f_2'(x)}{x^{\lambda_2-1}}$$
yields that $(f_1f_2)'(x)/x^{\lambda_!+\lambda_2-1}$ is again completely monotonic and this concludes the proof.\hfill $\square$

As mentioned in the introduction, a function $f$ belong to $\mathcal T_{\lambda,\alpha}$ if and only if $f'(x)/x^{\lambda-1}$ belongs to $\mathcal S_{\lambda+1-\alpha}$. We show that also $f(x)/x^{\lambda}$ is a generalized Stieltjes function of positive order:
\begin{prop}
    Let $f\in \mathcal T_{\lambda,\alpha}$. Then $(f(x)-f(0))/x^{\lambda}$ belongs to $\mathcal S_{\lambda+1-\alpha}$.
\end{prop}
{\it Proof.} Since $f'(x)/x^{\lambda-1}$ belongs to $\mathcal S_{\lambda+1-\alpha}$ we have, for some completely monotonic function $\varphi$, 
$$
f'(s)=cs^{\lambda-1}+s^{\lambda-1}\int_0^{\infty}t^{\lambda-\alpha}\varphi(t)e^{-st}\, dt.
$$
Integration of this relation between $0$ and $x$ yields, after changing the order of integration and a change of variable ($u=ts/x$),
\begin{align*}
f(x)-f(0)&=cx^{\lambda}/\lambda+\int_0^{\infty}\int_0^xe^{-st}s^{\lambda-1}\, ds\, t^{\lambda-\alpha}\varphi(t)\, dt\\
&=cx^{\lambda}/\lambda+x^{\lambda}\int_0^{\infty}\int_0^te^{-xu}u^{\lambda-1}\, du\, t^{-\alpha}\varphi(t)\, dt.
\end{align*}
This gives
$$
(f(x)-f(0))/x^\lambda=
c/\lambda+\int_0^{\infty}\sigma(u)u^{\lambda-\alpha}e^{-xu}\, du,
$$
where the function $\sigma$ is given by 
$$\sigma(u)=\int_1^{\infty}v^{-\alpha}\varphi(uv)\, dv.
$$
It is immediate that $\sigma$ is completely monotonic, and thus the proposition is proved.\hfill $\square$

As an application of this result we notice the following corollary.
\begin{cor}
\label{cor:product-higher-order}
    If $f_1\in \mathcal T_{\lambda_1,\alpha_1}$ and $f_2\in \mathcal T_{\lambda_2,\alpha_2}$ then $(f_1-f_1(0))(f_2-f_2(0))$   belongs to $\mathcal T_{\lambda_1+\lambda_2,\alpha_1+\alpha_2-1}$.
\end{cor}
{\it Proof.} We may assume that the functions are zero at the origin. Thus $f_j(x)/x^{\lambda_j}$ and $f_j'(x)/x^{\lambda_j-1}$ belong to $\mathcal S_{\lambda_j+1-\alpha_j}$. The elementary relation 
$$
\frac{(f_1f_2)'(x)}{x^{\lambda_!+\lambda_2-1}}=
\frac{f_1'(x)}{x^{\lambda_!-1}}\frac{f_2(x)}{x^{\lambda_2}}+
\frac{f_1(x)}{x^{\lambda_!}}\frac{f_2'(x)}{x^{\lambda_2-1}}$$
yields that $(f_1f_2)'(x)/x^{\lambda_!+\lambda_2-1}\in \mathcal S_{\lambda_1+1-\alpha_1}\mathcal S_{\lambda_2+1-\alpha_2}\subseteq \mathcal S_{\lambda_1+\lambda_2+1-(\alpha_1+\alpha_2-1)}$. This concludes the proof.\hfill $\square$

Let us mention an example related to the incomplete Beta function.
\begin{cor}
    Assume $0<b<a$. The function $B(b,a-b)-B(b,a-b,e^{-x})$ belongs to $\mathcal B_{[a-b]+1}$.
\end{cor}
{\it Proof.} 
Denote the function by $f$. By definition, 
$$
f(x)=\int_{e^{-x}}^1t^{b-1}(1-t)^{a-b-1}\, dt=\int_0^xe^{-bt}(1-e^{-t})^{a-b-1}\, dt,
$$
so that 
$$
x^{-[a-b]}f'(x)=e^{-bx}\left(\frac{1-e^{-x}}{x}\right)^{[a-b]}(1-e^{-x})^{a-b-[a-b]-1}.
$$
The function $1/(1-e^{-x})$ is logarithmically completely monotonic. Indeed,
$$
-\left(\log \frac{1}{1-e^{-x}}\right)'=\frac{e^{-x}}{1-e^{-x}}=\sum_{k=0}^{\infty}e^{-(k+1)x}.
$$
Therefore any positive power of it is completely monotonic. In particular $(1-e^{-x})^{a-b-[a-b]-1}$ is completely monotonic, and so $x^{-[a-b]}f'(x)$ is also completely monotonic. This completes the proof.
\hfill $\square$

\begin{rem}
 The function from the corollary above does not belong to any of the classes $\mathcal T_{[a-b]+1,\alpha}$, where $\alpha<[a-b]+2$. Indeed if it did belong to $\mathcal T_{[a-b]+1,\alpha}$ then the function $x^{-[a-b]}e^{-bx}(1-e^{-x})^{a-b-1}$ would be a generalized Stieltjes function of some positive order $\mu$, and consequently, 
$$
x^{\mu-[a-b]}e^{-bx}(1-e^{-x})^{a-b-1}
$$
would be increasing. This is a contradiction since the limit of the expression above as $x$ tends to infinity is $0$.   
\end{rem}

\begin{prop} 
    If $f\in \mathcal B_1$ and $f(x)/x$ is logarithmically completely monotonic then $f^{\lambda}\in \mathcal B_{\lambda}$ for any $\lambda\geq 1$.
\end{prop}
{\it Proof.} We have 
$x^{1-\lambda}\left(f(x)^\lambda\right)'=\lambda(f(x)/x)^{\lambda-1}f'(x)$,
which is completely monotonic since it is a product of completely monotonic functions.\hfill $\square$
\begin{prop}
    If $f\in \mathcal B_{\mu}$ for some $\mu\in (0,1]$ and $\lambda>0$ then $f^{\lambda}\in \mathcal B_{\mu([\lambda]+1)}$.
\end{prop}
{\it Proof.} We have
$$
\frac{\left(f(x)^\lambda\right)'}{x^{\mu([\lambda]+1)-1}}=\lambda\left(\frac{f(x)}{x^{\mu}}\right)^{[\lambda]}\frac{f'(x)}{x^{\mu-1}}\, f(x)^{\lambda-1-[\lambda]},
$$
which is a product of completely monotonic functions, since $f\in \mathcal B_1$ and the exponent $\lambda-1-[\lambda]$ is negative.\hfill $\square$
\begin{rem}
 Returning to the function $w$ we see that
$$
W(t)=t^{a-b-[a-b]-1}w(t)=e^{-bt}\left(\frac{1-e^{-t}}{t}\right)^{[a-b]}
(1-e^{-t})^{(a-b)-[a-b]-1}$$
is completely monotonic 
since $t\mapsto 1-e^{-t}$ is in $\mathcal B_1$. 
Since $\Gamma(x+b)/\Gamma(x+a)$ can be represented as 
 \begin{align*}
    \frac{\Gamma(x+b)}{\Gamma(x+a)}&=\frac{1}{\Gamma(a-b)}\int_0^{\infty}e^{-xt}t^{a-b-1}w(t)\, dt\\
    &=\frac{1}{\Gamma(a-b)}\int_0^{\infty}e^{-xt}t^{[a-b]}W(t)\, dt
\end{align*} this reproves that the ratio of Gamma functions belongs to $\mathcal S_{[a-b]+1}$.
\end{rem}

\begin{rem}
    The representation \eqref{eq:Mbdd} yields, with $\lambda=a-b$,
    $$
    \frac{L_{a,b}(x)}{x^{\lambda+1}}=\int_0^{\infty}p(t)t^{\lambda-1}e^{-xt}\, dt
    $$
     where $p(t)=\int_0^tu^3\varphi(u)\, du$, $\varphi$ being completely monotonic. It is clear that
     $$
     p(t)t^{-4}=\int_0^1s^3\varphi(ts)\, ds,
     $$
     and thus that $p$ is completely monotonic of order $-4$. Thus the function $L_{a,b}(x)/x^{a-b+1}$ belongs to $\mathcal S_{a-b+4}$.
\end{rem}

\begin{defn}
Let $\alpha\geq 0$. A function $f:(0,\infty)\to (0,\infty)$ belongs to the class $\mathcal C_\alpha$ if $x^\alpha f(x)$ is the Laplace transform of a decreasing and logarithmically convex function $u$. 
\end{defn}
We remark that for a positive and decreasing function $u$, $\mathcal L(u)$ is defined if and only if $\int_0^1u(t)\, dt<\infty$.

\begin{prop}
    \label{prop:stamatis-1}
    If $f\in \mathcal C_{\alpha}$ then 
    \begin{enumerate}[(a)]
        \item $x^{\alpha+1}f(x)\in \mathcal B_1$;
        \item $1/f(x)\in \mathcal B_{\alpha+1}$;
        \item $1/(x^{\alpha+1}fx))$ is completely monotonic.
    \end{enumerate}
\end{prop}
\noindent {\it Proof.} By definition $x^{\alpha}f(x)=\mathcal L(u)(x)$, where $u$ is decreasing and logarithmically convex. Since $u$ is decreasing, it follows that $x^{\alpha+1}f(x)=x\mathcal L(u)(x)$ is a Bernstein function (see e.g.\ the remark below). This proves (a). According to \cite[Theorem 11.11]{S}, $x\mathcal L(u)(x)$ is actually a so-called special Bernstein function. Therefore also the function $1/\mathcal L(u)$ is a Bernstein function.  From Proposition \ref{prop:product-bernstein} it now follows that
$$
\frac{1}{f(x)}=x^{\alpha}\, \frac{1}{\mathcal L(u)(x)}\in \mathcal B_{\alpha+1},$$
and this proves (b). Item (c) follows from (b) and \cite[Corollary 2.1]{KP1}, and this concludes the proof. \hfill $\square$
\begin{rem}
    Suppose that $u(t)=\int_t^{\infty}d\mu(s)$. Then by Fubini's theorem, 
    $$
    x\mathcal L(u)(x)= \int_0^{\infty}\left(1-e^{-xs}\right)\, d\mu(s),
    $$
    showing that $x\mathcal L(u)(x)$ is a Bernstein function.

\end{rem}

We end this paper by applying Proposition \ref{prop:stamatis-1} to the $\psi$-function. 

\begin{cor}
    For $0<b<a-1$ the function $F_{a,b}$ defined as
    $$
    F_{a,b}(x)=\frac{x^{a-b-1}}{\psi(x+a)-\psi(x+b)}
    $$
    belongs to $\mathcal B_{a-b}$.
\end{cor}
\noindent {\it Proof.} 
 Notice that $\psi(x+a)-\psi(x+b)=\mathcal L(\varphi_{a,b})(x)$, where $\varphi_{a,b}$ is the decreasing and logarithmically convex function in the lemma below. Hence, $\mathcal L(\varphi_{a,b})(x)/x^{a-b-1}$ belongs to $\mathcal C_{a-b-1}$, and thus by Proposition \ref{prop:stamatis-1}, its reciprocal $x^{a-b-1}/\mathcal L(\varphi_{a,b})(x)$ belongs to $\mathcal B_{a-b}$. The proof is complete.\hfill $\square$

\begin{lemma}
    \label{lemma:logconv}
    For $0<b<a-1$ the function $\varphi_{a,b}(t)=(e^{-bt}-e^{-at})/(1-e^{-t})$ is decreasing and logarithmically convex.
\end{lemma}
\noindent {\it Proof.} Since $\varphi_{a,b}(t)=\Phi_{a,b}(t)-(b-a)$ it is clear from Lemma \ref{lemma:Phi} and Lemma \ref{lem:eta} that it decreases. The logarithmic convexity can be seen as follows. Letting $\lambda=a-b$, a standard computation shows that logarithmic convexity holds if and only if 
$$
\frac{e^{-t}}{(1-e^{-t})^2}>\frac{\lambda^2 e^{-\lambda t}}{(1-e^{-\lambda t})^2}.
$$
Rearranging the terms this is equivalent to 
\begin{equation}
\label{eq:ineq}
\left(\frac{\sinh (\lambda t/2)}{\sinh (t/2)}\right)^2>\lambda^2.
\end{equation}
However, by comparing coefficients in the Taylor series it easily follows that $\sinh (\lambda t/2)/\sinh (t/2)>\lambda$, and therefore \eqref{eq:ineq} holds.\hfill $\square$

\noindent
Stamatis Koumandos\\
Department of Mathematics and Statistics\\
The University of Cyprus\\
P. O. Box 20537\\
1678 Nicosia, Cyprus\\
email: skoumand@ucy.ac.cy
\medskip

\vspace{0.1in}

\noindent
Henrik Laurberg Pedersen\\
Department of Mathematical Sciences\\
University of Copenhagen\\
Universitetsparken 5\\
DK-2100, Denmark\\
email: henrikp@math.ku.dk
    
\end{document}